\documentclass[11pt]{article}
\usepackage{amssymb}
\usepackage{latexsym}
\newtheorem{theorem}{Theorem}[section]
\newtheorem{prop}{Proposition}[section]

\newtheorem{example}{Example}[section]

\newtheorem{lemma}{Lemma}[section]

\newcommand{\qed}{\hspace*{\fill}\hbox{$\Box$}\vspace{2ex}}

\newcommand{\be}{\begin{equation}}
\newcommand{\ee}{\end{equation}}
\newcommand{\bea}{\begin{eqnarray}}
\newcommand{\eea}{\end{eqnarray}}
\newcommand{\bean}{\begin{eqnarray*}}
\newcommand{\eean}{\end{eqnarray*}}
\newcommand{\benu}{\begin{enumerate}}
\newcommand{\eenu}{\end{enumerate}}
\newcommand{\edo}{\end{document}}

\newcommand{\barr}{\begin{array}}
\newcommand{\earr}{\end{array}}

\newcommand{\CC}{{\mathbb C}}

\newcommand{\HH}{{\mathbb H}}

\newcommand{\RR}{{\mathbb R}}

\newcommand{\Ric}{{\rm Ric}}

\renewcommand{\Re}{{\rm Re}}
\newcommand{\SL}{{\rm SL}}

\newcommand{\mcr}{\mathcal{R}}
\newcommand{\ecke}{\;_-\!\rule{0.2mm}{0.2cm}\;\;}

\begin{document}

\title{The twistor equation in Lorentzian spin geometry}
\author{Helga Baum and Felipe Leitner}

\date{}
\maketitle

\begin{abstract}
In this paper we discuss the twistor equation in Lorentzian spin geometry.
In particular, we explain the local conformal structure of Lorentzian
manifolds, which admit twistor spinors inducing lightlike Dirac currents.
Furthermore, we derive all local geometries with singularity free  
twistor spinors that occur up to dimension 7.\\


\end{abstract}

\section{Introduction}

A classical object in differential geometry are conformal Killing fields. 
These are by definition infinitesimal
conformal symmetries, i.e. the flow of such vector fields preserves the conformal 
class of the metric. The number of
linearly independent conformal Killing fields measures the degree of conformal 
symmetry on the manifold. This number
is bounded by $\frac{1}{2}(n+1)(n+2)$, where $n$ is the dimension of the manifold. 
If it is the maximal one the manifold is conformally flat. 
S. Tachibana and T. Kashiwada (cf. \cite{Tachibana/Kashiwada:69}, \cite{Kashiwada:68})
introduced a generalization of conformal Killing fields, the conformal Killing forms 
(or twistor forms). Conformal Killing forms are solutions of a conformally invariant
 twistor type equation on differential forms. They were studied
in General Relativity mainly from the local viewpoint in order to
integrate the
 equation of motion (e.g. \cite{Penrose/Walker:70}), furthermore they were used to 
obtain symmetries of field equations (\cite{Benn/ua1:97},
\cite{Benn/ua2:97}).  Recently, U. Semmelmann (\cite{Semmelmann:01}) started to 
discuss global properties of conformal Killing forms in Riemannian geometry. 
Another generalization of conformal Killing vectors is that of conformal Killing
spinors (or twistor spinors), which are solutions of the conformally invariant twistor
 equation on spinors introduced by R. Penrose in General Relativity 
(cf. \cite{Penrose/Rindler:86}).  Whereas conformal Killing fields are classical
symmetries, conformal Killing spinors define infinitesimal symmetries on supermanifolds 
(cf. \cite{Alekseevski/Cortes/ua:98}). Special kinds of such spinors, parallel and special 
Killing spinors, occur in supergravity and string theories. 
In 1989 A. Lichnerowicz and Th.  Friedrich started a systematic study of twistor
spinors in conformal Riemannian geometry. Whereas the global structure of Riemannian 
manifolds admitting 
twistor spinors is quite well understood (cf. e.g. \cite{Lichnerowicz1:88}, 
\cite{Lichnerowicz2:88}, \cite{Lichnerowicz1:89}, \cite{Friedrich:89},
 \cite{Lichnerowicz2:90}, \cite{Baum/Friedrich/ua:91},
\cite{Habermann:90}, \cite{Habermann:93}, \cite{Habermann:94}, \cite{Habermann:96}, 
\cite{Kuehnel/Rademacher:94}, \cite{Kuehnel/Rademacher1:96}, \cite{Kuehnel/Rademacher1:97}, 
\cite{Kuehnel/Rademacher:98}), the state of art in its origin, Lorentzian geometry, 
is far from being satisfactory. We are mainly interested in the following problems:
\begin{enumerate} \item {\em Which Lorentzian geometries admit twistor spinors?} 
\item {\em How are the properties of twistor spinors related to the geometric structures 
where they can occur?} \end{enumerate} J. Lewandowski (\cite{Lewandowski:91}) described 
the local normal forms of 4-dimensional spacetimes with zero free twistor spinors.
His results indicated that there are interesting relations between 
twistor spinors, different global contact structures and Lorentzian geometry that should 
be discovered. H. Baum (\cite{Baum:99}, \cite{Baum:00}) described twistor spinors on Fefferman 
spaces and on Lorentzian symmetric spaces. Ch. Bohle (\cite{Bohle:03}) and F. Leitner 
(\cite{Leitner:03}) studied Lorentzian geometries with Killing spinors, which are a special class 
of twistor spinors. In the present paper we describe the geometric 
structures, which appear up to dimension 7.  We consider only the case of zero free 
twistor spinors. Results for twistor spinors with zeros can be found in \cite{Leitner:01}.\\
 After recalling the definition of twistor spinors we
discuss in section 3 Brinkmann spaces, Lorentzian Einstein--Sasaki
structures and Fefferman spaces and their relation to the problem in question. 
In chapter 4 we study the local conformal structure of 
Lorentzian manifolds that admit twistor spinors inducing lightlike Dirac currents 
(cf. Proposition \ref{prop-9}). 
In chapter 5 we derive all local conformal structures of Lorentzian
manifolds admitting singularity free solutions of the twistor equation in 
low dimensions $n \leq 7$ (cf. Theorem \ref{theo}).\\

\section{The twistor equation on spinors}

In this section we recall the definition of twistor spinors and fix some notations. 
For more details we refer to \cite{Penrose/Rindler:86} or \cite{Baum/Friedrich/ua:91}.\\
Let $(M^n,g)$ be a semi-Riemannian spin manifold of dimension $n \ge 3$. We denote by $S$
 the spinor bundle and by $\mu :T^*M \otimes S \to S$ the
Clifford multiplication. The 1-forms with values in the spinor bundle
 decompose into two subbundles
\begin{displaymath}
T^* M \otimes S = V \oplus Tw,
\end{displaymath}
where $V$, being the orthogonal complement to the `twistor bundle' $Tw := Ker \mu\;$, 
is isomorphic to $S$.
Usually, we identify $TM$
and $T^*M$ using the metric $g$. \\
We obtain two differential operators of first order by composing the spinor
derivative $\nabla^S$ with the orthogonal projections onto each of these
subbundles,\\
the {\em Dirac operator} $D$
\begin{displaymath}
D: \Gamma (S) \stackrel{\nabla^S}{\longrightarrow} \Gamma (T^* M \otimes S)=
\Gamma (S \oplus Tw) \stackrel{pr_S}{\longrightarrow}\Gamma (S)
\end{displaymath}
and the {\em twistor operator} $P$
\begin{displaymath}
P: \Gamma (S) \stackrel{\nabla^S}{\longrightarrow} \Gamma (T^*M \otimes S) =
\Gamma (S \oplus Tw) \stackrel{pr_{Tw}}{\longrightarrow} \Gamma (Tw).
\end{displaymath}
Locally, these operators are given by the following formulas
\begin{eqnarray*}  
D \varphi &=& \sum\limits^n_{i=1} \sigma^i \cdot \nabla^S_{s_i} \varphi\\
P \varphi &=& \sum\limits^n_{i=1} \sigma^i \otimes (\nabla_{s_i}^S \varphi +
\frac{1}{n} s_i \cdot D \varphi),
\end{eqnarray*}
where  $(s_1 , \ldots , s_n)$ is a local orthonormal basis, $(\sigma^1,
\ldots, \sigma^n)$ its dual and $\cdot$ denotes the Clifford multiplication.
 Both operators are conformally
covariant. More exactly, if $\tilde{g} = e^{2 \sigma } g$ is a conformal
change of the metric, the Dirac and the twistor operator satisfy
\begin{eqnarray*}  
D_{\tilde{g}} &=& e^{- \frac{n+1}{2} \sigma} D_g e^{\frac{n-1}{2} \sigma}\\ 
P_{\tilde{g}} &=& e^{- \frac{\sigma}{2}} P_g e^{- \frac{\sigma}{2}}.
\end{eqnarray*}  
A spinor field is called {\em twistor spinor}  or {\em conformal Killing spinor} if it 
lies in the kernel of the
twistor operator $P$.
Using the local formula for the twistor operator one obtains the following characterization 
of a twistor spinor:
A spinor field $\varphi \in \Gamma(S)$ is a twistor spinor if and only if
\begin{displaymath}
\nabla^S_X \varphi + \frac{1}{n} X \cdot D\,\varphi = 0  \qquad  \mbox{for all vector fields}
\;\; X.
\end{displaymath}
Obviously, each parallel spinor ($\nabla^S\varphi=0$) is a twistor spinor.
An other special class of twistor spinors are the Killing spinors $\varphi$,
which satisfy $\nabla_X^S \varphi = \lambda X\cdot\varphi$ for some $\,\lambda\in\CC \setminus \{0\}$. 
It is a well-known fact, that -- as in the case  of conformal vector fields -- 
the dimension of the space of twistor
spinors is bounded and the maximal possible  dimension is attained only for conformally 
flat manifolds. More exactly,
it holds\\[-0.2cm]
\begin{prop} (cf. \cite{Baum/Friedrich/ua:91})\label{prop-2}
\begin{enumerate}
\item The dimension of the space of twistor spinors is a conformal
invariant and bounded by
\begin{displaymath}
\dim \ker P \le 2 \,\mathrm{rank} S = 2^{[\frac{n}{2}]+1} =: d_n.
\end{displaymath}  
\item If  $\dim \ker P =d_n$ then
$(M^n,g)$ is conformally flat.
\item If $(M^n,g)$ is simply connected and conformally flat then $\dim \ker
P= d_n$.
\end{enumerate}
\end{prop}
\ \\
Hence, for example, all simply connected space forms $\RR^n_k, \tilde{\HH}^n_k
, \tilde{S}^n_k$ admit the maximal number of linearly independent twistor
spinors.

\section{Twistor spinors on Lorentzian spin manifolds - Examples}

Now, we restrict our attention to the case of Lorentzian signature $(-+
\ldots +)$. In this section we will explain the special geometries that occur in 
the Theorem \ref{theo}.\
Let $(M^n,g)$ be an oriented and time-oriented Lorentzian spin
manifold. On the spinor bundle $S$ there exists an indefinite non-degenerate
inner product $\langle \cdot ,\cdot \rangle$ such that
\begin{eqnarray*}
\langle X \cdot \varphi  , \psi \rangle & =& \langle \varphi , X \cdot \psi
\rangle \quad \quad \mbox{and}\\
X(\langle \varphi , \psi \rangle )&= &\langle \nabla^S_X \varphi , \psi
\rangle + \langle \varphi , \nabla^S_X \psi \rangle,
\end{eqnarray*}
for all vector fields $X$ and all spinor fields $\varphi, \psi$ 
(cf. \cite{Baum:81}). Each
spinor field $\varphi \in \Gamma (S)$ defines a vector field $V_{\varphi}$
on $M$, the so-called Dirac current, by
\begin{equation} \label{gl03}
g (V_{\varphi} , X) := - \langle X \cdot \varphi , \varphi \rangle .
\end{equation}
A direct calculation shows the following properties of the Dirac current\\[-0.2cm]

\begin{prop} \label{prop-3} (cf. \cite{Baum:99})
Let $\varphi$ be a spinor field on a Lorentzian spin manifold $(M^n,g)$
with Dirac current $V_{\varphi}$. Then
\begin{enumerate}
\item $V_{\varphi}$ is causal and future-directed.
\item The zero sets of $\varphi$ and $V_{\varphi}$ coincide.
\item If $\varphi$ is a twistor spinor, $V_{\varphi}$ is a conformal Killing
field.
\end{enumerate}
\end{prop}
\ \\
Now, let us discuss 3 types of special Lorentzian geometries that admit
twistor spinors.

\subsection{Brinkmann spaces with parallel spinors}

A Lorentzian manifold is called \emph{Brinkmann space} if it admits a
non-trivial  lightlike parallel vector field. Let us consider two examples 
of such spaces.\\[-0.3cm]
\begin{example}[pp-manifolds.]
{\em A Brinkmann space is called  \emph{pp-manifold} if its  
Riemannian curvature tensor $\mcr$ satisfies
\begin{displaymath}
\mbox{Trace}\,_{(3,5),(4,6)} \mcr \otimes \mcr =0 .
\end{displaymath}
Equivalently, pp-manifolds can be characterized as those Lorent\-zian
mani\-folds $(M^n,g)$, where the metric has the following local normal form
depending only on one function $f$ of $(n-1)$ variables
\begin{displaymath}
g= dt\, ds +f(s, x_1 , \ldots , x_{n-2}) ds^2 + \sum\limits^{n-2}_{i=1}
dx_i^2 .
\end{displaymath}
(cf. \cite{Schimming:74}).
In terms of holonomy, pp-manifolds can be characterized as those
Lorentzian manifolds, for which the restricted holonomy group is contained in the 
Abelian normal subgroup $1 \ltimes \RR^{n-2}$ of the parabolic subgroup 
$(\RR \times SO(n-2)) \ltimes \RR^{n-2}$ in  $SO(1,n-1)$
(see \cite{Leistner:01}, \cite{Leistner:02}).
Using the latter fact one can easily prove that for each simply connected
pp-manifold
\begin{displaymath}
\dim \ker P \ge \frac{d_n}{4} .
\end{displaymath}
Furthermore, on generic pp-manifolds each twistor spinor is parallel. An
important example of geodesically complete pp-manifolds are the Lorentzian
symmetric spaces with solvable transvection group (Cahen-Wallach-spaces) 
(cf. \cite{Cahen/Wallach:70}, \cite{Neukirchner:02}).
}
\end{example}
\ \\[-0.8cm]
\begin{example}[Brinkmann spaces with special K\"ahler flag]
{\em Let $(M^n,g)$ be a \\
Brink\-mann space with the lightlike parallel vector field
$V$. Then $V$ defines a flag of subbundles $\,\RR V \subset V^{\perp} \subset TM \,$ in $TM$,
where $V^{\perp} =\{ Y \in TM \,|\, g(V,Y)=0 \}$. We equip the bundle $E:= V^{\perp} /\RR V$
with the
positive definite inner product $\tilde{g}$ induced by $g$ and with the  metric
connection $\tilde{\nabla}$ induced by the Levi-Civita connection of $g$. We
call $\RR V \subset V^{\perp} \subset TM$ a \emph{
K\"ahler flag}, if in case of even $n$ there is a parallel orthogonal almost complex structure
 $J: E \longrightarrow E$  and if in case of odd $n$ 
there exists a parallel subbundle $H \subset E$ of codimension 1 equipped
with a parallel orthogonal almost complex stucture $J: H \longrightarrow H$. \\
The K\"ahler flag $\; \RR V \subset V^{\perp} \subset TM\;$ is called {\em special K\"ahler flag} if in addition  
\[ \mbox{Trace}\; (J \circ \mcr^{\tilde{\nabla}} (X,Y)) =0\, \mbox{
for all}\;\, X,Y \in TM. \]
A Brinkmann space has a special K\"ahler flag iff its reduced holonomy representation 
is contained in $SU(2m-2) \ltimes \RR^{2m-2}\,$ (if $n=2m$) resp. in $\,(SU(2m-2)\times 1) 
\ltimes \RR^{2m-1}$ (if $n=2m+1$).
It was proved by I.Kath in \cite{Kath:99} that a Brinkmann space $(M^{n},g)$ has a special K\"ahler flag 
if and only if $(M^{n},g)$ has
pure parallel spinors.
}
\end{example}
\ \\[-0.8cm]
\subsection{Twistor spinors on Lorentzian Einstein--Sasaki manifolds}  

An odd-dimensional Lorentzian manifold $(M^{2m+1}, g ; \xi)$ equipped
with a vector field $\xi$ is called \emph{Lorentzian Sasaki
manifold} if
\begin{enumerate}
\item $\xi$ is a timelike Killing field with $g( \xi, \xi )=-1.$
\item The map $J:= - \nabla \xi :TM \to TM$ satisfies
\begin{eqnarray*}
J^2 X &=& -X -g (X, \xi ) \xi \quad \quad \mbox{and}\\
(\nabla_X J)(Y) &=& - g(X,Y) \xi + g (Y, \xi )X.
\end{eqnarray*}
\end{enumerate}
Let us consider the metric cone $\,C_- (M) := (\RR^+ \times M , -dt^2 + t^2 g)\,$ 
with timelike cone axis over
$(M,g)$. The cone metric has signature $(2,2m)$. Then the following relations
between properties of $M$ and those of its cone are easy to verify \\[-0.3cm]

\begin{tabbing}
xxxxxxxxxxxxxxxxxxxxxxxxxxxxxxxx\= xxxxxxxx\= \kill
$(M^{2m+1},g;\xi)$ \> \>  cone $C_-(M)$ \\[0.2cm]
 Lorentzian Sasaki \> $\Longleftrightarrow$ \> (pseudo)-K\"ahler\\[0.05cm]
 Lorentzian Einstein--Sasaki ($R<0$) \> $\Longleftrightarrow$ \>  
Ricci-flat and (pseudo)-K\"ahler\\[0.05cm]
 Lorentzian Einstein--Sasaki ($R<0$) \> $\Longleftrightarrow$\>  $\mathrm{Hol}_0 (C_- (M)) 
\subset
SU(1,m)$\\[-0.2cm]
\end{tabbing}
The standard example for regular Lorentzian Einstein--Sasaki manifolds are $S^1$-
bundles over Riemannian K\"ahler--Einstein spaces of negative scalar
curvature:
Let $(X^{2m},h)$ be a Riemannian K\"ahler--Einstein spin manifold of
scalar
curvature $R_X <0$ and let $(M^{2m+1}, \pi , X ; S^1)$ denote the
$S^1$-principal bundle associated to the square root $\sqrt{\Lambda^{m,0} X}$
of the canonical bundle of $X$ given by the spin structure. Furthermore,
let $A$ be the connection on $M$ induced by the Levi-Civita
connection of $(X,h)$. Then
\begin{displaymath}
g := \pi^* h - \frac{16m}{(m+1) R_X} A \odot A
\end{displaymath}
defines a Lorentzian Einstein--Sasaki metric on the spin manifold
$M^{2m+1}$. Lorent\-zian Einstein--Sasaki manifolds admit a special kind
of twistor spinors.\\[-0.2cm]

\begin{prop} (cf. \cite{Kath:99}, \cite{Bohle:03}) \label{prop-4}
Let $(M,g)$ be a simply connected Lorent\-zian Einstein--Sasaki manifold. Then
$M$ is spin and admits a twistor spinor $\varphi$ on $M$ such that
\begin{itemize}
\item[a)] $V_{\varphi}$ is a timelike Killing field with
$g(V_{\varphi} , V_{\varphi})=-1$
\item[b)] $V_{\varphi}  \cdot \varphi = - \varphi$
\item[c)] $\nabla_{V_{\varphi}} \varphi = - \frac{1}{2} i \varphi$.
\end{itemize}
Conversely, if $(M,g)$ is a Lorentzian spin manifold with a twistor spinor
satisfying the conditions a), b) and c). Then $\xi := V_{\varphi}$ is a Lorent\-zian
Einstein--Sasaki structure on $(M,g)$. Each twistor spinor $\varphi$
on a Lorentzian Einstein-Sasaki manifold is the sum of two Killing spinors.
\end{prop}
\ \\[-0.8cm]
\subsection{Twistor spinors on Fefferman spaces}

Fefferman spaces are Lorent\-zian manifolds which appear in the frame work
of CR geometry. Let us first explain the necessary notations from CR
geometry. Let $N^{2m+1}$ be a smooth oriented manifold of odd dimension
$2m+1$. A \emph{CR structure on $N$} is a pair $(H,J)$ where\\[0.2cm]
\hspace*{0.1cm} 1.  $H \subset TM$ is a real $2m$-dimensional subbundle.\\[0.1cm]
\hspace*{0.1cm} 2.  $J:H \to H$ is an almost complex structure on $H$, i.e. $J^2 =-Id$.\\[0.1cm]
\hspace*{0.1cm} 3.  If $X,Y \in \Gamma (H)$, then $[JX,Y]+[X,JY] \in \Gamma (H)$ and\\[0.1cm]
    \hspace*{1cm}  $J([JX,Y]+[X,JY]) - [JX, JY]+[X,Y] \equiv 0$  (integrability 
condition).\\[0.2cm]
Let $(N,H,J)$ be a CR manifold. In order to define Fefferman spaces we fix a contact form
$\theta \in \Omega^1 (N)$ on $N$ such that
$\theta|_H = 0$. Let us denote by $T$ the Reeb vector field of $\theta$,
which is defined by the conditions $\theta (T) =1$ and  $T \ecke d \theta =0$. In the
following we suppose that the Levi form $L_{\theta} :H \times H \to \RR$
\begin{displaymath}
L_{\theta} (X,Y) := d \theta (X,JY)
\end{displaymath}
is positive definite. Then $(N,H, J , \theta)$ is called a
{\em strictly pseudoconvex manifold}. The tensor field $g_{\theta} := L_{\theta}
+ \theta \odot \theta$ defines a Riemannian metric on $N$. There is a special
 metric covariant derivative on a strictly pseudoconvex manifold, the
{\em Tanaka-Webster connection} $\nabla^W : \Gamma (TN) \to \Gamma (TN^*
\otimes TN)$, uniquely defined by the conditions
\begin{eqnarray*}
\nabla^W g_{\theta} &=& 0\\
Tor^W (X,Y) &=& L_{\theta} (JX,Y) \cdot T\\
Tor^W (T,X) &=& - \frac{1}{2} ([T,X]+J[T,JX])
\end{eqnarray*}  
for $X,Y \in \Gamma (H)$. This connection satisfies $\nabla^W J= 0$ and
$\nabla^W T=0$ (cf. \cite{Tanaka:75}, \cite{Webster:78}). Let us denote by
$T_{10} \subset TN^{\CC}$ the eigenspace of the complex extension of $J$ on $H^{\CC}$
to the eigenvalue $i$. Then $L_{\theta}$ extends to a Hermitian form on $T_{10}$ by
\[L_{\theta} (U,V) := - i\,d\theta (U, \bar{V}), \;\;\;U,V \in T_{10}.\]
For a complex 2-form $\omega \in \Lambda^2 N^{\mathbb{C}}$ we denote by
$\mbox{trace}_{\theta} \omega$ the $\theta$-trace of $\omega$:
\begin{displaymath}
\mbox{trace}_{\theta} \omega := \sum\limits^m_{\alpha =1} \omega (Z_{\alpha} ,
\bar{Z}_{\alpha} ),
\end{displaymath}
where $(Z_1 , \ldots , Z_m)$ is a unitary basis of $(T_{10}, L_{\theta})$.
Let $\mathfrak{R}^W$ be the $(4,0)$-curvature tensor of the Tanaka-Webster   
connection $\nabla^W$ on the complexified tangent bundle of $N$
\begin{displaymath}
\mathfrak{R}^W (X,Y,Z,V) := g_{\theta} (([\nabla^W_X , \nabla^W_Y]-
\nabla^W_{[X,Y]}) Z, \bar{V}) .
\end{displaymath}
and let us denote by
\begin{displaymath}
\Ric^W := \mbox{trace}^{(3,4)}_{\theta} := \sum\limits^m_{\alpha =1} \mathfrak{R}^W
(\cdot , \cdot , Z_{\alpha} , \bar{Z}_{\alpha})
\end{displaymath}
the \emph{Tanaka-Webster Ricci curvature} and by $R^W := \mbox{trace}_{\theta}
\Ric^W$ the \emph{Tanaka-Webster scalar curvature}. The Ricci curvature $\Ric^W$ is a
$(1,1)$-form on $N$ with $\Ric^W (X,Y) \in i \RR$ for real vectors $X,Y \in
TN$. The scalar curvature $R^W$ is a real function.\\
Now, let us suppose that $(N^{2m+1}, H, J, \theta)$ is a strictly
pseudoconvex spin manifold. The spin structure of $(N, g_{\theta})$ defines   
a square root $\sqrt{\Lambda^{m+1,0} N}$ of the canonical line bundle
\begin{displaymath}
\Lambda^{m+1,0} N:= \{ \omega \in \Lambda^{m+1} N^C \,|\, V \ecke
\omega =0 \,\,\,\forall V \in \bar{T}_{10} \} .
\end{displaymath}
We denote by $(F, \pi , N)$ the $S^1$-principal bundle associated to
$\sqrt{\Lambda^{m+1,0} N}$.
Let $A^W$ denote the connection form on $F$ defined by the Tanaka-Webster
connection $\nabla^W$. Then 
\begin{displaymath} 
h_{\theta} := \pi^* L_{\theta} - i \,\frac{8}{m+2} \pi^* \theta \odot
(A^W-\frac{i}{4(m+1)}R^W\cdot\theta)\ 
\end{displaymath}
is a Lorentzian metric such that the conformal class $[h_{\theta}]$ is an
invariant of the CR structure $(N,H,J)$. The metric $h_{\theta}$ is
$S^1$-invariant, the fibres
of the $S^1$-bundle are lightlike. We call
$(F^{2m+2}, h_{\theta})$ with its canonically induced spin structure
\emph{Fefferman space of the strictly pseudoconvex spin manifold $(N,H,J,
\theta)$}.\\[-0.2cm]

\begin{prop} (\cite{Baum:99}) \label{prop-5}
Let $(N,H,J, \theta)$ be a strictly pseudoconvex spin manifold with the
Fefferman space $(F, h_{\theta})$. Then there exist two linearly independent
twistor spinors $\varphi$ on $(F, h_{\theta})$ such that
\noindent
\begin{itemize}
\item[a)] $V_{\varphi}$ is a regular lightlike Killing field
\item[b)] $V_{\varphi} \cdot \varphi =0$
\item[c)] $\nabla_{V_{\varphi}} \varphi = i\, c\, \varphi$, where
$c \in \RR \backslash \{ 0 \}$.
\end{itemize}
Conversely, if $(M,g)$ is an even dimensional Lorentzian spin
manifold with a twistor spinor satisfying a), b) and c), then there
exists a strictly pseudoconvex spin manifold $(N,H,J, \theta)$ such
that its Fefferman space is locally isometric to $(M,g)$.
\end{prop}
\ \\[-0.5cm]
\section{Twistor spinors inducing lightlike Dirac currents}

As we noticed in Proposition \ref{prop-3}, each twistor spinor $\varphi$
induces a causal conformal vector field $V_{\varphi}$. 
In this section we study the case that the Dirac current $V_\varphi$ is lightlike. 
Aiming at a local conformal classification we may assume in addition that $V_\varphi$ is Killing.\\
Let us start with some notations.
Let $W : \Lambda^2 M \to \Lambda^2 M\,$ denote the Weyl tensor of $(M^n,g)$ considered
as selfadjoint map on the space
of 2-forms and
$\,\mathrm{Ric}$ denotes the Ricci tensor of $(M^n,g)$ considered here as
$(1,1)$-tensor or as $(2,0)$-tensor whatever is needed.
In conformal geometry, there are two further curvature tensors that play an
important role, the Rho tensor $K$
\begin{displaymath}
K(X) := \frac{1}{n-2} \Big( \frac{R}{2(n-1)} X - \Ric (X) \Big)\ ,\quad
X \in TM
\end{displaymath}
and the Cotton-York tensor
\begin{displaymath}
C(X,Y):= (\nabla_X K)(Y) - (\nabla_Y K)(X) , \quad X,Y \in TM.
\end{displaymath}
For twistor spinors we have the following properties of the Cotton-York and 
the Weyl tensor\\[-0.2cm]
\begin{prop} \label{prop-6}
Let $\varphi$ be an arbitrary twistor spinor. Then the Dirac current $V_{\varphi}$
annihilates the Cotton-York and the Weyl tensor:
$\; V_{\varphi} \ecke C = 0$,\ $V_{\varphi} \ecke W = 0$.
\end{prop}
\begin{proof} We use the following well-known integrability conditions for
twistor spinors $\varphi$ (cf. \cite{Baum/Friedrich/ua:91})
\begin{eqnarray}
\label{Int3} W (\eta )
\cdot \varphi &=& 0 \quad \mbox{for all 2-forms}\;\; \eta\\
\label{Int4}
W(X \wedge Y)\cdot D\varphi& = & n\, C(X,Y) \cdot \varphi\,.
\end{eqnarray}
We deduce from this
\begin{eqnarray*}
C(V_{\varphi},X,Y) &=& g(V_{\varphi},C(X,Y)) = - \langle C(X,Y) \cdot \varphi,\varphi 
\rangle \\
                   &=& -\frac{1}{n} \langle W(X \wedge Y) \cdot  D\varphi, \varphi 
\rangle 
= \frac{1}{n}
\langle D\varphi, W(X\wedge Y) \cdot \varphi \rangle = 0.
\end{eqnarray*}
Moreover, with the relation $X\cdot\eta=-X\ecke\eta+X^\flat\wedge\eta$ in the Clifford algebra,
where $X$ denotes a vector and $\eta$ a $2$-form, we have 
\begin{eqnarray*}
W(V_\varphi,X,Y,Z)=-\langle\varphi,W(X,Y,Z)\cdot\varphi\rangle=-\langle\varphi,
Z^\flat\wedge W(X,Y)\cdot\varphi\rangle\in\RR\ .
\end{eqnarray*}
Since $\langle\varphi,\rho^3\cdot\varphi\rangle\in i\RR$ for all $3$-forms $\rho^3$,
it follows that $V_\varphi \ecke W=0$.
\qed
\end{proof}
\ \\[0.2cm]
Now, let us mention a special property of spinor fields with lightlike Dirac current:
\begin{prop} (cf. \cite{Leitner:01}) \label{lem-1}
Let $\varphi$ be a spinor field on a Lorentzian manifold with lightlike Dirac current 
$V_{\varphi}$. Then
\begin{enumerate}
\item
$\; V_{\varphi} \cdot \varphi = 0.$
\item
$\; \langle \varphi, \varphi \rangle = 0.$
\end{enumerate}
\end{prop}
\begin{proof} 
The claimed properties for spinors in the lemma are purely algebraic.
Therefore, it is sufficient and appropriate to prove these properties on the level of
the corresponding representations. For this
we use the usual concrete realization of the representation
of the Clifford algebra
$Cl_{1,n-1}$ and its complexification $Cl_{1,n-1}^{\CC}$ on the spinor
module $\Delta_{1,n-1}$
in terms of Kronecker products of matrices (\cite{Baum/Friedrich/ua:91}). 
Let us consider the complex $(2\times 2)$-matrices
\[
E := \left(\small \begin{array}{cc}1 & 0\\0 & 1 \end{array}\right),\qquad T := 
\left(\small
\begin{array}{cc} 0 & -i\\ i & 0 \end{array}\right), \qquad g_1 :=
\left(\small\begin{array}{cc}i
& 0\\ 0 & -i\end{array}\right),\qquad g_2 := \left(
\small\begin{array}{cc}0 & i\\ i & 0\end{array}\right)
\]
and let $\tau(1)=i$ and $\tau(2)=\tau(3)=\dots =\tau(n)=1$.
We denote by $(e_1,\ldots,e_n)$ an orthonormal basis of the 
Minkowski space $\RR^{1,n-1}$.
If $n = 2m$ then a Clifford representation on the
spinor module $\Delta_{1,n-1}\cong  \CC^{2^{m}}$
is realized by the map
$
\Phi_{2m} : Cl^{\CC}_{1,n-1} \to \CC (2^m)
$, which is 
generated by
\[
\begin{array}{rcl}
\Phi_{2m}(e_{2j-1}) & = & \tau (2j-1) \cdot E \otimes\cdots\otimes E
\otimes g_1 \otimes \underbrace{T\otimes\cdots\otimes
T}_{(j-1)-times}\\[8mm]
\Phi_{2m}(e_{2j}) & = & \tau (2j)\cdot E \otimes\cdots\otimes E
\otimes g_2
\otimes \underbrace{T\otimes\cdots\otimes T}_{(j-1)-times}\ ,
\end{array}
\]
where $j = 1,\cdots,m$. If $n=2m + 1$ then a representation
$ 
\Phi_{2m+1}: Cl^{\CC}_{1,n-1} \to \CC (2^m) 
$
is generated by 
\[
\begin{array}{rcl}
\Phi_{2m+1}(e_j) & = & \Phi_{2m}(e_j)\qquad \mbox{for}\ j=1,\ldots, 2m\quad \mbox{and}\\[2mm]
\Phi_{2m+1}(e_n) & = & i\cdot \;T \otimes\cdots\otimes T.
\end{array}
\]
Furthermore, let us consider the vectors 
$\,u(\nu)=\frac{1}{\sqrt{2}}\left(\scriptsize\begin{array}{c}1\\-i 
\nu\end{array}\right)\in\CC^{2},\, \nu=\pm 1,\,$
 and the unitary basis 
$\,
\{\; u(\nu_1,\ldots,\nu_m):=u(\nu_1)\otimes\cdots\otimes u(\nu_m)\mid\ \nu_i\in\{\pm 1\}
\ \}
\,$
of $\Delta_{1,n-1}\cong\CC^{2^m}$ with respect to the standard scalar product
$(\cdot,\cdot)_\Delta$ of $\CC^{2^m}$. The indefinte scalar product
$\langle\cdot,\cdot\rangle$ on $S$ is defined by the $Spin_o(1,n-1)$-invariant
inner product $\langle v,w\rangle_\Delta=(e_1\cdot v,w)_{\Delta}$ on $\Delta_{1,n-1}$.
Let 
$\ell:\Delta_{1,n-1}\rightarrow\RR^{1,n-1}$ denote the 
map, which maps a 
spinor
$v$ to its Dirac current $\ell(v)=\langle v,e_1v\rangle e_1-\sum_{i\geq 2}
\langle v,e_iv\rangle e_i$. 
We calculate now the inverse image 
$\ell^{-1}(\RR(e_1+e_2))$ of the lightlike direction $\RR(e_1+e_2)$. 
For this let
\[v=\sum_{(\nu_1,\ldots,\nu_m)\in\{\pm 1\}^m}   
a_{\nu_1,\ldots,\nu_m}\cdot
u(\nu_1,\ldots,\nu_m),\qquad
a_{\nu_1,\ldots,\nu_m}\in\CC,\]
be an arbitrary spinor represented in the unitary basis of
$\Delta_{1,n-1}$. It holds
\[\begin{array}{rcl}
e_1\cdot v&=&-\sum
a_{\nu_1,\ldots,\nu_m}\cdot
u(\nu_1,\ldots,-\nu_m)\qquad\mbox{and}\\[1mm]
e_2\cdot v&=&\ \ \sum
\nu_m \cdot a_{\nu_1,\ldots,\nu_m}
u(\nu_1,\ldots,-\nu_m).\end{array}\]
Then we obtain
\[\begin{array}{l}\langle v,e_1v\rangle_{\Delta}=(e_1\cdot v,e_1\cdot v)_\Delta=
\sum |a_{\nu_1,\ldots,\nu_m}|^2\qquad\mbox{and}\\[1.5mm]
\langle v,e_2v\rangle_{\Delta}=(e_1\cdot v,e_2\cdot v)_\Delta=
\sum -\nu_m\cdot|a_{\nu_1,\ldots,\nu_m}|^2\ .\end{array} \]
It is $\ell(v)\in\RR(e_1+e_2)$ if and only if
$\langle v,e_1v\rangle=-\langle v,e_2v\rangle$. The latter condition
is equivalent to
$\,
a_{\nu_1,\ldots,\nu_{m-1},-1}=0\,$ 
for all 
$\,(\nu_1,\ldots,\nu_{m-1})\in\{\pm 1\}^{m-1}$.
Hence, a spinor $v$ with $\ell(v)\in\RR(e_1+e_2)$
has the form
\[
v=a\otimes u(1), \qquad a\in\bigotimes_{m-1}\CC^2.\]
Then $\langle v,v\rangle_\Delta=(e_1\cdot v,v)_\Delta=
-(a\otimes u(-1), a\otimes u(1))=0$ and $(e_1+e_2)\cdot v=a\otimes (ig_1+g_2)
u(1)=0$.
This proves the desired properties in case that $\ell(v)\in\RR(e_1+e_2)$.
Since the map $\ell$ is equivariant under the action of the spin group $Spin_o(1,n-1)$ 
and 
the spin group acts transitive on the lightlike directions
in the lightcone of the Minkowski space $\RR^{1,n-1}$, we can conclude 
that the claimed properties for spinors 
$v\in\Delta_{1,n-1}$ 
with 
arbitrary lightlike Dirac current 
are true in general. 
\qed
\end{proof}
\ \\[2mm]
Let $V$ be a vector field and let $\theta$ denote the dual 1-form $\theta(X)=g(X,V)$.  
Then $V$ is called {\em
twisting}
if the 3-form $d\theta \wedge \theta$  nowhere vanishes and {\em non-twisting} if 
$d\theta \wedge \theta=0$. By the
Frobenius Theorem the latter means  that the distribution $V^{\perp} \subset TM$ is 
integrable. 
For twistor spinors with non-twisting Dirac current we have the following result.\\[-0.2cm]
\begin{prop}(cf. \cite{Leitner:01}) \label{prop-8}
Let $\varphi$ be a twistor spinor with  lightlike, non-twisting Dirac
current. Then $(M,g)$ is locally conformally
equivalent to a Brinkmann space with parallel spinor.
\end{prop}
\begin{proof} 
The condition that $V_{\varphi}$ 
has no zero and no twist
implies by the Frobenius Theo\-rem that locally there
are functions $\sigma,f$ such that
$V_{\varphi}=e^{-2\sigma}\mathrm{grad}f$. Then, locally we have 
$V_{\varphi}=\widetilde{\mathrm{grad}}f$
with respect to the metric $\tilde{g}=e^{2\sigma}g$. Hence, without loss of
generality, we may
assume that $V:=V_{\varphi}$ is a lightlike conformal gradient field. Let
$V=\mathrm{grad}f$. Then
\[ \mathrm{Hess}f = \frac{1}{2}L_Vg = \frac{\Delta f}{2n}g \qquad
\mathrm{and} \qquad 
\nabla_X V = \frac{\mathrm{div}V}{2n}X \;\; \mbox{for all} \, X.
\]
 Since $V$ is lightlike, we obtain
$0= X(g(V,V))=2 g(\nabla_XV,V)=\frac{\mathrm{\mathrm{div}V}}{n}g(X,V)$ for all $X$ and
 therefore div$(V)=0$.
 This shows that 
the Dirac current $V_{\varphi}$ is parallel. Now we show that the spinor
$\varphi$ is parallel, too.
To this end, let $(s_1,\ldots,s_n)$ be a local orthonormal frame with
$V_{\varphi}=s_1+s_2$. By Proposition \ref{lem-1} we know that
$V_{\varphi}\cdot \varphi = 0$. Since $V_{\varphi}$ 
is parallel, we obtain from $\,0= \nabla_X(V_{\varphi}\cdot \varphi) =
\nabla_X V_{\varphi} \cdot \varphi 
+ V_{\varphi}\cdot \nabla_X \varphi\,$ that
\[ V_{\varphi} \cdot \nabla_X \varphi =0 \quad \mbox{and}\mbox \quad s_1
\cdot s_2 \cdot \nabla_X \varphi =
 - \nabla_X \varphi \quad \mbox{for all}\; X. \]  
Let $X$ be a vector field with $g(X,X)=\pm 1$ and let $\psi:= g(X,X) X \cdot
\nabla_X \varphi\,$. 
Since $\varphi$ is a twistor spinor, the spinor $\psi$ does not depend on the choice of $X$. 
Choose $X \in V_{\varphi}^{\perp}$. 
Then $V_{\varphi}\cdot \psi = -g(X,X)X \cdot V_{\varphi}\cdot \nabla_X \varphi =0.\,$ 
On the other hand, for $X=s_1$
we obtain
\[ 0= V_{\varphi} \cdot \psi = -(s_1+s_2) \cdot s_1 \cdot \nabla_{s_1} \varphi = 
(-1+s_1\cdot s_2)\cdot \nabla_{s_1} \varphi = -2 \nabla_{s_1}\varphi. \]
This shows that $\psi=0$, which implies $\nabla \varphi =0$.  
\qed
\end{proof}
\ \\[0.2cm]
Using the conditions on the twist one is able to  characterize -- at least locally -- 
all geometries that admit a
zero free
twistor spinor with lightlike Dirac  current. In particular, the
following Proposition explains the role, Fefferman
spaces and Brinkmann spaces are playing among all Lorentzian geometries that admit 
twistor spinors.\\
\begin{prop} \label{prop-9}
Let $(M,g)$ be a Lorentzian spin manifold admitting a twistor spinor
$\varphi$ such that $V_{\varphi}$ is lightlike and Killing. Then the function
$\Ric (V_{\varphi} , V_{\varphi})$ is constant and non-negative on $M$.
Furthermore,
\begin{enumerate}
\item $\Ric (V_{\varphi}, V_{\varphi}) >0$ if and only if  
$(M,g)$ is
locally isometric to a Fefferman space. In this case the Dirac current is twisting and 
the dimension of $M$ is even.
\item $\Ric (V_{\varphi} , V_{\varphi}) =0$ if and only if $(M,g)$ is
locally conformal equivalent to a Brinkmann space with parallel spinors. In this case 
the Dirac current is
non-twisting .
\end{enumerate}
\end{prop}
\begin{proof} At first, we observe the following general relations
for a lightlike Killing field $V$.
It is $g(\nabla_X\nabla_YV,Z)=\mcr(Y,Z,X,V)$ for all
vector fields $X,Y$ and $Z$ in $TM$, especially
$g(\nabla_XV,\nabla_YV)=\mcr(V,X,Y,V)$. If in addition $V \ecke W=0$
then \begin{eqnarray*} g(\nabla_X\nabla_YV,Z) &=& g\star K(Y,Z,X,V)\\
&=& -g(Y,V)K(Z,X)-g(Z,X)K(Y,V)\\
& & +g(Y,X)K(Z,V)+g(Z,V)K(Y,X)\ ,
\end{eqnarray*}
where $\star$ denotes the Kulkarni-Nomizu product.
Next we show that
\[
X(Ric(V,V))= 4 K(\nabla_XV,V)=0\qquad \mbox{for\ all}\quad
X\in TM\ .\] 
Fix $x \in M$ and let $(e_1,\ldots,e_n)$ denote an orthonormal frame, arising by
parallel displacement from $x$. Then in the point $x$ we have the following identities. 
Using the skew-symmetry of $\,g(\cdot,\nabla V)\,$ and the symmetry of $K$ we
obtain
\begin{eqnarray}
\sum_kg(e_k,e_k)\cdot
K(\nabla_{e_k}V,e_k) = 0 .\label{Gr1}
\end{eqnarray} 
The second Bianchi identity for the Riemannian curvature tensor
and $V(R)=0$
yields 
\begin{eqnarray}
\sum_kg(e_k,e_k)\,e_k\left(
K(V,e_k)\right) = 
\sum_k g(e_k,e_k)\cdot (\nabla_{e_k}K)(V,e_k) = 0 .
\label{Gr2}
\end{eqnarray}
 Furthermore, from 
 $V\ecke C=0$ follows
\begin{eqnarray}
V(K(e_i,V))=
(\nabla_V K)(e_i,V)= (\nabla_{e_i}K)(V,V) \ .\label{Gr3}
\end{eqnarray} 
It is
\begin{eqnarray*}
\mcr(e_i,\nabla_{e_k}V,e_k,V)+\mcr(e_k,\nabla_{e_i}V,e_k,V)
&=&g(\nabla_{e_k}\nabla_{e_i}V,\nabla_{e_k}V)
+g(\nabla_{e_k}\nabla_{e_k}V,\nabla_{e_i}V)\\ &=&e_k(
\mcr(V,e_i,e_k,V))\ , \end{eqnarray*}
\begin{eqnarray*}
\mcr(e_i,\nabla_{\!e_k}\!V,e_k,V)\!+
\!\mcr(e_k,\nabla_{\!e_i}\!V,e_k,V)\!\!\!
&=\ \
&\!\!\!\!\!g(e_i,e_k)K(\nabla_{e_k}V,V)-g(e_i,V)K(\nabla_{e_k}V,e_k)\\
&\ \
+&\!\!\!\!\!g(e_k,e_k)K(\nabla_{e_i}V,V)-g(e_k,V)K(\nabla_{e_i}V,e_k)\\
&\ \ -&\!\!\!\!\! g(\nabla_{e_i}V,e_k)K(e_{k},V)\qquad\qquad\quad 
\mbox{and}\end{eqnarray*}
\begin{eqnarray*}e_k(
\mcr(V,e_i,e_k,V))&=\ \ &\!\!\!\!\!
g(V,e_k)\cdot e_k(K(e_i,V))+g(e_i,\nabla_{e_k}V)\cdot K(V,e_k)\\
&\ \ +&\!\!\!\!\!g(e_i,V)\cdot e_k(K(V,e_k))
-g(e_i,e_k)\cdot e_k(K(V,V))\ .\end{eqnarray*}
Summing up the latter equations and using (\ref{Gr1}), (\ref{Gr2}) and
(\ref{Gr3}) results in
\begin{eqnarray*}
\sum_k g(e_k,e_k)\left(
\mcr(e_i,\nabla_{e_k}V,e_k,V)+\mcr(e_k,\nabla_{e_i}V,e_k,V)\right)
&=&(n-1) K(\nabla_{e_i}V,V)\qquad \mbox{and}\\
\sum_k g(e_k,e_k)\cdot   
e_k(
\mcr(V,e_i,e_k,V)) &=& -3K(\nabla_{e_i}V,V)\ .
\end{eqnarray*}
Hence, $\; K(\nabla_{e_i} V,V) = 0 $.
Eventually,
\begin{eqnarray*}e_i(Ric(V,V))&=&2\sum_kg(e_k,e_k)\cdot
g(\nabla_{e_i}\nabla_{e_k}V,\nabla_{e_k}V)\\&=&
2\sum_kg(e_k,e_k) g\star K(e_k,\nabla_{e_k}V,e_i,V)\\
&=&4 K(\nabla_{e_i}V,V)\ .\end{eqnarray*}
In particular, $Ric(V,V)$ is constant on $M$.\\
Now, let us consider $V=V_{\varphi}$. Then the condition   
$\Ric(V,V)=\mathrm{const}\;$ 
follows from Proposition \ref{prop-6}. We denote by
$\eta$ and $\theta$ the 1-forms
\[
\eta(X):=K(V,X) ,\qquad \theta(X):= g(V,X).
\]
Furthermore, let $T$ be the vector field dual  to $\eta$. For the following, we 
normalize the spinor $\varphi$
such
that
$\eta(V)= -(n-2)\,\Ric(V,V)=: \varepsilon \in \{0,1,-1\}$. Let us consider the 
endomorphism
\[
J:TM \longrightarrow TM \qquad J(X):= \nabla_X V\ .
\]
The map $J$ is skew-adjoint and satisfies $JV=JT=0$, since $V_\varphi$ is
lightlike Killing 
and $g(JT,X)=-K(V,\nabla_XV)=0$ for all $X$.
Then we obtain for arbitrary vectors $X$ and $Y$ 
\begin{eqnarray*}
g(J^2(X),Y) &=& -g(J(X),J(Y))=-\mcr(V,X,Y,V)\\
&=& -g(V,Y)K(X,V)-g(X,V)K(V,Y)\\
 & & 
 +g(V,V)K(X,Y)+g(X,Y)K(V,V)\\
&=& g(X,y)\eta(V) - g(V,Y)\eta(X) -\theta(X)g(T,Y),
\end{eqnarray*}
which shows that
\begin{eqnarray}
J^2(X) &=& \varepsilon X - \theta(X)T - \eta(X)V \label{J1}\ .
\end{eqnarray}
Moreover, it is
\begin{eqnarray}
g(T,V) &=& \varepsilon\;,\qquad g(V,V)\;=\;g(T,T)\;=\;0\label{J2}\\
d\theta(X,Y)&=& 2 g(JX,Y) \label{J3}
\end{eqnarray}
Now, let $\varepsilon=1$. We consider the $J$-invariant subbundle 
$\,H=\mathrm{span}\{V_{\varphi},T)^{\perp} \subset TM\,$. 
The spin structure of $(M,g)$ reduces to a spin structure of the hermitian bundle 
$(H,J,g)$. The spin structure, $H$, $\theta$ and $J$ project down to the (locally defined) 
manifold $N$ resulting from $M$ by factoring out the integral curves of $V$. Then it can 
be proved as in  
\cite{Graham:87} that $(N,H,J,\theta)$ is a strictly pseudoconvex manifold and that 
 $(M,g)$ is 
locally
isometric to the Fefferman
space of $(N,H,J,\theta)$. The case $\varepsilon=-1$ can not 
occur,  since by
(\ref{J1}) and (\ref{J2}) $J$ would be a skew-adjoint
involution on the positive definite subbundle $\,\mathrm{span}(V,T)^{\perp}\,$ of $TM$.
In case that $\varepsilon=0$ the vector field $T$
is parallel to $V$. This implies that $Im(J^2)\subset \RR V$. Since $Im J \subset V^{\perp}$ 
we obtain $J(V^{\perp}) \subset \RR V$. Then (\ref{J3}) shows that
 $d\theta$ vanishes on $V^\bot$, which means that
$\theta$ is
non-twisting.  Hence, the second assertion of 
the Theorem follows from Proposition \ref{prop-8}. 
\qed
\end{proof}
\section{Twistor spinors in dimension $n \le 7$}
In this section we discuss the twistor equation on Lorentzian spin manifolds in the 
dimensions $n=3,4,5,6$ and
$7$. For solutions of the twistor equation without singularities we derive
a complete list of possible underlying local Lorentzian geometries in these low
dimensions.  For dimension $n=4$ this was
already
proved by J. Lewandowski. To begin with, we state some special properties
for spinors on low dimensional Lorentzian
spin manifolds. These properties can be derived from the representation
theory of the spinor modules in low dimensions 
(cf. \cite{Lawson/Michelsohn:89}) and the discussion
of their orbit structure, which can be found e.g. in \cite{Bryant:00}.

\begin{lemma}\label{lem-2} {\em(Properties of spinors on low dimensional
Lorentzian manifolds)}
\begin{enumerate}
\item Let $n=3$. Then there exists a real structure $\tau: S
\longrightarrow S$ on the spinor
 bundle such that
$\, \tau(X \cdot \varphi)= - X \cdot \tau(\varphi)\,$  and  $\; \langle
\tau \varphi,\tau \psi \rangle = -\langle
\psi,\varphi \rangle.$ Each spinor field $\varphi \in \Gamma(S)$
satisfies
 $V_{\varphi} \cdot \varphi = \langle \varphi, \varphi \rangle \varphi
\,$ and 
 $\,g(V_{\varphi},V_{\varphi})=-\langle \varphi,\varphi
\rangle^2.$
If $\varphi$
  is a real spinor field,  i.e. $\tau(\varphi)=\varphi$,
then the Dirac current $V_{\varphi}$ is lightlike.
\item  Let $n=5$. Then there exists a quaternionic structure
$J: S \longrightarrow S$ such that 
 $ \,J(X \cdot \varphi) = X \cdot J(\varphi)\,$ and  $\,
\langle J\varphi,J\psi \rangle = -\langle \psi, \varphi
\rangle.$ Again each spinor field $\varphi \in \Gamma(S)$
satisfies
 $V_{\varphi} \cdot \varphi = \langle \varphi, \varphi \rangle \varphi
\,$ and 
 $\,g(V_{\varphi},V_{\varphi})=-\langle \varphi,\varphi
\rangle^2.$
\item Let $n=7$. Then there exists
a quaternionic structure
$J: S \longrightarrow S$ such that   
 $ \,J(X \cdot \varphi) = -X \cdot J(\varphi)\,$ and
 $\,
\langle J\varphi,J\psi \rangle =  \langle \psi, \varphi
\rangle.$ Each spinor field $\varphi \in \Gamma(S)$
satisfies
$V_{\varphi}\cdot\varphi=\langle\varphi,\varphi\rangle\varphi 
+\langle\varphi,J\varphi\rangle J\varphi$ and
$|\langle\varphi,\varphi\rangle|^2+|\langle\varphi,J\varphi\rangle|^2=
- g(V_{\varphi},V_\varphi)$.
\item  Let $n=2,4,6$. Then each spinor field $\varphi \in \Gamma(S^{\pm})$
satisfies
$\,V_{\varphi} \cdot \varphi =0\,$ and  $\,
g(V_{\varphi},V_{\varphi})=0 $. 
\end{enumerate}
\end{lemma}

\begin{proof}The existence
of the real or quaternionic structures is clear from
the representation theory of the spinor modules
(cf. \cite{Lawson/Michelsohn:89}). 
The main point here is to prove the formulas
concerning
$V_{\varphi}\cdot\varphi$ and $g(V_{\varphi},V_{\varphi})$. Although these
formulas are very
natural, it seems that there is no natural proof for them.\\ 
Therefore, 
we 
explain in the following the orbit structure of the spinor modules
$\Delta_{1,n-1}$ for $n=3,4,5,6$ and $7$ with respect to the action
of the spin group and calculate explicitly
with respect to 
convenient normal forms of representatives
in the various orbits. 
Let $n=3$. There is a real structure $\tau$ on 
$\Delta_{1,2}\cong\CC^2$,
which is invariant under the action of $Spin_o(1,2)$ and
$Cl_{2,1}\subset Cl^{\CC}_{1,2}$. 
The real spinor representation
$\Delta_{1,2}^{\RR}$
is isomorphic to the standard representation of $\SL(2,\RR)$ on 
$\RR^2$. Beside the zero orbit, there are the orbit types
to the representatives
\[
\sigma_1=
\left(\scriptsize\begin{array}{c}1\\id\end{array}
\right),\qquad
\sigma_2=
\left(\scriptsize\begin{array}{c}1+ic\\0\end{array}\right)\ 
\qquad\mbox{and}\qquad
\sigma_3=
\left(\scriptsize\begin{array}{c}i\\0\end{array}\right)\ 
\]
parametrized by $0\neq c\in\RR$ and $d\in\RR$.
We choose the Clifford representation generated by 
\[
\Phi(e_1):=\left(\scriptsize
\begin{array}{cc}0&i\\-i&0\end{array}\right),\qquad
\Phi(e_2):=\left(\scriptsize\begin{array}{cc}0&-i\\-i&0
\end{array}\right)\qquad
\mbox{and}\qquad
\Phi(e_3):=\left(\scriptsize
\begin{array}{cc}i&0\\0&-i\end{array}\right)\ .\]
Then it is
\[V_{\sigma_1}
=\left(\scriptsize\begin{array}{c}\ \ 1+d^2\\-1+d^2\\0
\end{array}\right),\qquad
V_{\sigma_2}
=\left(\scriptsize\begin{array}{c}\ \ 1+c^2\\-1-c^2\\0\end{array}\right)
\qquad\mbox{and}\qquad
V_{\sigma_3}
=\left(\scriptsize\begin{array}{r}1\\-1\\0\end{array}\right)\in\RR^{1,2}
\]
and calculating the Clifford product gives 
$V_{\sigma_i}\cdot\sigma_i=\langle\sigma_i,
\sigma_i\rangle_\Delta\cdot\sigma_i$\
for $i=1,2,3$. In particular, 
$g(V_{\sigma_i},V_{\sigma_i})=-\langle\sigma_i,
\sigma_i\rangle^2$ and 
every real spinor has lightlike
Dirac current.\\ 
Let $n=5$. There is a $Cl_{1,4}$-invariant quaternionic
structure on the spinor module
$\Delta_{1,4}\cong\HH^2$. The 
orbit types are determined
by the representatives 
\[
\sigma_1=
\left(\scriptsize\begin{array}{c}r\\0\end{array}
\right),\qquad
\sigma_2=
\left(\scriptsize\begin{array}{c}0\\r\end{array}\right)\
\qquad\mbox{and}\qquad
\sigma_3=
\left(\scriptsize\begin{array}{c}1\\1\end{array}\right)\in\HH^2,
\]
where $r\in\RR_+$. 
We generate the Clifford representation by
\[\begin{array}{l}
\Phi(e_1):=\left(\scriptsize
\begin{array}{cc}1&0\\0&-1\end{array}\right),\qquad
\Phi(e_2):=\left(\scriptsize\begin{array}{cc}0&1\\-1&0
\end{array}\right),\\[4mm]\Phi(e_3):=\left(\scriptsize
\begin{array}{cc}0&i\\i&0\end{array}\right),\qquad
\Phi(e_4):=\left(\scriptsize
\begin{array}{cc}0&
j\\j&0\end{array}\right)\qquad\mbox{and}\qquad
\Phi(e_5):=\left(\scriptsize
\begin{array}{cc}0&
k\\k&0\end{array}\right),\end{array}\]
where $i,j,k$ denote the imaginary units.
Then
$V_{\sigma_1}=V_{\sigma_2}=r^2\cdot e_1$
and $V_{\sigma_3}=e_1+e_2$. Executing the Clifford multiplication
results in
$V_{\sigma_i}\cdot\sigma_i=\langle\sigma_i,
\sigma_i\rangle_\Delta\cdot\sigma_i$, which also shows that
$g(V_{\sigma_i},V_{\sigma_i})=-\langle\sigma_i,
\sigma_i\rangle^2$.\\
In general, the even-dimensional spinor representations $\Delta_{1,2m-1}$
split into the half spinor representations $\Delta_{1,2m-1}^{\pm}$.
The half spinor modules are spanned by
\[\Delta_{1,2m-1}^{\pm}=\{u(\nu_1,\ldots,\nu_m):\
\Pi_{i=1}^m\nu_i=\pm 1\} \  . \]
In particular, this shows that in all even dimensions there are  positive
and negative half spinors inducing
lightlike Dirac currents
(cf. Proof of Proposition
\ref{lem-1}).\\  
For $n=2$ the half spinors are represented by
$\varphi=r\cdot u(1)$ and $r\cdot u(-1)$, where $r\in\RR_+$. These spinors
have
lightlike
Dirac current and Proposition \ref{lem-1} implies that $V_\varphi\cdot
\varphi=0$.\\
For $n=4$ the real half spinor representations
$\Delta^{\pm}_{1,3}$
are isomorpic to the canonical representation of $\SL(2,\CC)$ 
on $\CC^2\cong\RR^4$. Hence, there is exactly one non-trivial orbit in
each
half spinor module. These two orbits are represented by
$u(1,1)$ and $u(1,-1)$, which give rise to lightlike Dirac currents.\\
Let $n=6$. There is a $Cl_{5,1}$-invariant quaternionic structure on
$\Delta_{1,5}\cong\HH^4$. The occuring orbit types are
\[\sigma_1=\left(\scriptsize\begin{array}{c}1\\0\\0\\0\end{array}\right),\qquad
\sigma_2=\left(\scriptsize\begin{array}{c}0\\0\\1\\0\end{array}\right),\qquad
\sigma_3=\left(\scriptsize\begin{array}{c}1\\r\\0\\0\end{array}\right)\qquad\mbox{and}
\qquad
\sigma_4=
\left(\scriptsize\begin{array}{c}1\\0\\ \lambda\\0\end{array}\right),\]
where $r\in\RR_+$ and $\lambda\in\HH$. Both half spinor modules admit
exactly one non-trivial orbit and these are represented
by $\sigma_1$ and $\sigma_2$, which then must have lightlike Dirac
current.\\
Finally, let $n=7$. The module  $\Delta_{1,6}$ admits a
$Cl_{6,1}$-invariant
quaternionic structure $J$ and the
restriction of $\Delta_{1,6}$ to the action
of $Spin(1,5)$ is isomorphic to the spinor representation $\Delta_{1,5}$.
>From 
the orbit type classification of $\Delta_{1,5}$ we derive the 
orbit types
in $\Delta_{1,6}$, which are then parametrized by
\[\sigma_\lambda=
\left(\scriptsize\begin{array}{c}1\\0\\\lambda\\0\end{array}\right),\quad
\lambda=\lambda_1+\lambda_2j\in Im\HH\ .\] 
It is
$\langle\sigma_\lambda,\sigma_\lambda\rangle=-2i\lambda_1$. Moreover,
$J\sigma_\lambda=
\left(\scriptsize\begin{array}{c}j\\0\\-\lambda_2+\lambda_1j\\
0\end{array}\right)$ and
$\langle\sigma_\lambda,J\sigma_\lambda\rangle=2i\overline{\lambda_2}$.
For a suitable realization of the Clifford representation
we have $V_{\sigma_\lambda}=(1+|\lambda|^2)e_1+(1-|\lambda|^2)e_2$,
where $e_1$ is a timelike unit vector,
and
$V_{\sigma_\lambda}\cdot\sigma_\lambda=\sigma_\lambda\cdot(-2i\lambda)$.
This proves that
$V_{\sigma_\lambda}\cdot\sigma_\lambda=
\sigma_\lambda\cdot\langle\sigma_\lambda,\sigma_\lambda\rangle+
J\sigma_\lambda\cdot\langle\sigma_\lambda,J\sigma_\lambda\rangle$.
\qed
\end{proof}

\begin{lemma} \label{lem-3}
The dimension of the space of twistor spinors on a non-conformally flat Lorentzian
 manifold $(M^n,g)$ is bounded
by\\[0.1cm]
(1) $\;\;n=3$: $\qquad\dim_{\CC}\;\mathrm{Ker}\,P \leq 1$.\\
(2) $\;\;n=4$: $\qquad\dim_{\CC}\;\mathrm{Ker}\,P \leq 2$.\\
(3) $\;\;n=5$: $\qquad\dim_{\CC}\;\mathrm{Ker}\,P \leq 2$.
\end{lemma}
\begin{proof} Using the twistor equation one shows that in case $\,\mathrm{dim}_{\CC}\mathrm{Ker}P > 1\;$ 
there exists a dense set $A \subset M$ such that 
$\;\mathrm{dim}_{\CC}\{ \varphi(x) \in S_x \ | \ \varphi \in \mathrm{Ker}P \} \geq 2\;$ for all $x \in A$.\\
In dimension 3, the Weyl tensor $W$ vanishes. Hence, for each twistor
spinor $\varphi$ and all vector fields  
$X,Y$ the condition
$\; C(X,Y) \cdot \varphi = 0 \;$ 
holds (cf. (\ref{Int4})). Let $\mathrm{dim}_{\CC}\mathrm{Ker}P > 1\;$ and 
let $x \in A$ be an arbitrary point. Since $\mathrm{dim}_{\CC}S_x=2$, 
the vectors 
$C_x(U,W)$ annihilate $S_x$ 
for all $U, W \in T_xA$. This implies $C=0\;$.
A 4-dimensional Lorentzian manifold with a twistor spinor $\varphi$ is of Petrov type N or 0 in each point $x \in M$ 
where $\varphi(x) \not = 0$. Let $\varphi^+$ be the positive part of $\varphi$ and suppose $\varphi^+(x) \not = 0$.
Then $W_x=0$ or $L_x:= \RR V_{\varphi^+}(x)$ is the uniquely determined 4-fold principal null direction of $W_x$ and 
$L_x \cdot \varphi^+(x) = 0\;$. Let $\,\mathrm{dim}_{\CC}\mathrm{Ker}P > 2\;$. Then without loss of generality 
we may assume that $\,\mathrm{dim}_{\CC}(\mathrm{Ker}P^+:=\mathrm{Ker}P \cap \Gamma(S^+)) > 1\;$. Consider a dense 
set $A \subset M$ with   
$\;\mathrm{dim}_{\CC}\{ \varphi^+(x) \in S^+_x \ | \ \varphi^+ \in \mathrm{Ker}P^+ \} = \mathrm{dim}_{\CC}S^+_x = 2\;$.
Assume $x \in A$ and $W_x\not = 0$. Then the 4-fold principle null direction $L_x$ of $W_x$ annihilates $S_x^+$.
Hence $W=0$.
In dimension 5, $\mathrm{Ker}P$ is a quaternionic space (Lemma \ref{lem-2}). If 
$\,\mathrm{dim}_{\CC}\mathrm{Ker}P > 2\;$, than $\,\mathrm{dim}_{\HH}\mathrm{Ker}P > 1\;$ and there exists a dense 
set $A \subset M$ such that
$\;\mathrm{dim}_{\HH}\{ \varphi(x) \in S_x \ | \ \varphi \in \mathrm{Ker}P \} = \mathrm{dim}_{\HH}S_x = 2\;$.  
The integrability condition $\; W(\eta) \cdot \varphi = 0$ for all 2-forms $\eta$ shows that the 2-forms $W_x(\eta_x)$ 
annihilate $S_x$ for all $x \in A$. Hence, $W=0$.
\qed   
\end{proof}

We call a twistor spinor $\varphi$  {\em singularity free} if it has no zeros and the Dirac current $V_{\varphi}$ 
does not changes the causal type.  
There are the following geometric
structures
of Lorentzian manifolds with singularity free twistor spinors in dimension $n \leq 7$.
\\[-0.0cm]
\begin{theorem} (\cite{Lewandowski:91}, \cite{Leitner:01}) \label{theo}
Let $(M^n,g)$ be a Lorentzian manifold with a singularity free twistor spinor.
Then $(M^n,g)$ is locally conformally equivalent to one of the following
kinds of Lorentzian structures.\\[0.2cm]
%
\begin{tabular}{ll}
$n=3:$ & $\bullet$ pp-manifold \\[-0.2cm] & \\
$n=4:$ & $\bullet$ pp-manifold\\
& $\bullet$ Fefferman space \hspace{7cm}\\[-0.2cm] & \\
$n=5:$ & $\bullet$ pp-manifold\\
& $\bullet$ Lorentzian Einstein--Sasaki manifold\\
& $\bullet$ $\RR^{1,0}\times (N^4,h)$, where $(N^4,h)$ is Riemannian
Ricci-flat
K\"ahler\\[-0.2cm] & \\
$n=6:$ & $\bullet$ pp-manifold\\
& $\bullet$ Fefferman space\\
& $\bullet$ $\RR^{1,1} \times (N^4,h)$, where $(N^4,h)$ is Riemannian 
Ricci-flat K\"ahler\\
& $\bullet$ Brinkmann space with special K\"ahler flag \\[-0.2cm] &\\
$n=7:$ &  $\bullet$ pp-manifold \\
&  $\bullet$ Lorentzian Einstein-Sasaki manifold\\  
&  $\bullet$ $\RR^{1,2} \times (N^4,h)$, where $(N^4,h)$ is Riemannian Ricci-flat K\"ahler\\
&  $\bullet$ $\RR^{1,0} \times (N^6,h)$, where $(N^6,h)$ is Riemannian Ricci-flat K\"ahler\\
&  $\bullet$ $\RR^{1,0} \times (N^6,h)$, where $(N^6,h)$ is nearly K\"ahler, non-K\"ahler\\ 
&  $\bullet$ Brinkmann space with special K\"ahler flag \\ 
\end{tabular}
\end{theorem}
\begin{proof}
Let $n=3$ and let denote by $\tau: S \longrightarrow S$ the real structure on $S$ (Lemma \ref{lem-2}).
 If $\varphi$ is a zero  free twistor spinor then locally the real part 
$\Re\varphi :=
\frac{1}{2}(\varphi + \tau(\varphi))$  or the imaginary part $\Im \varphi := 
\frac{i}{2}(\tau(\varphi)-\varphi)$ is a
zero free real spinor. Hence,
according to Lemma \ref{lem-2}, we
have a zero free twistor spinor with
lightlike Dirac current, which becomes a Killing vector field after a conformal change of the
metric. Since the dimension is odd, $(M,g)$  is locally conformally equivalent to a Brinkmann space 
with parallel spinor
(Proposition
\ref{prop-9}). The 3-dimensional Brinkmann spaces with parallel spinors
have restricted holonomy group contained in 
$\,1 \ltimes \RR \subset SO(1,2)$, hence they are pp-manifolds (cf. Example 3.1). \\
In dimension $n=4$ and $n=6$ we may  assume that $\varphi$ is a zero free
half-spinor in $S^{\pm}$, since we are
only
interested in local considerations.
Since $V_{\varphi}$ has no zero (cf. Proposition \ref{prop-3}(2)) 
and is lightlike
(cf. Lemma \ref{lem-2}) there is a  local conformal transformation of the metric 
such that $V_{\varphi}$ is Killing
and
lightlike with respect to the conformally changed  metric.
Hence, according to Proposition \ref{prop-9}, $(M,g)$ is locally
conformal
equivalent to a Fefferman space or to a  Brinkmann space with parallel spinors. \\
The restricted holonomy group of a 4-dimensional Brinkmann space with parallel spinor 
is contained in $\,1\ltimes \RR^2 \subset SO(1,3)$, hence it is 
a pp-manifold. \\
The restricted holonomy group of a 6-dimensional Brinkmann space with parallel spinor is contained in 
$1 \ltimes \RR^4 \subset SO(1,5)$, $SU(2)\ltimes \RR^4 \subset SO(1,5)$ or 
in $1 \times SU(2)\subset SO(1,1) \times SO(4) \subset SO(1,5)$. Hence it is a pp-manifold, 
a Brinkmann space with special K\"ahler flag or the metric product of the hyperbolic plane $\RR^{1,1}$
with a 4-dimensional Riemannian Ricci-flat K\"ahler manifold. \\
Now, let us consider $n=5$. Let $\varphi$ be a zero free twistor spinor 
such that $V_{\varphi}$ is a lightlike or a timelike vector field. 
In case that
the Dirac current $V_{\varphi}$ is lightlike we may assume as above that $V_{\varphi}$ is 
Killing as well. Then, since the dimension is odd, 
Proposition \ref{prop-9} shows  that $(M,g)$ is locally conformal equivalent to 
a Brinkmann space with parallel spinors. The restricted holonomy group of a 5-dimensional Brinkmann spaces 
with parallel spinors
is contained in $1 \ltimes \RR^3 \subset SO(1,4)$, hence it is a pp-manifold.
Now, suppose that $V_{\varphi}$ 
is timelike. Then, by Lemma \ref{lem-2}, the length 
function $\,\langle \varphi,\varphi \rangle\,$ of $\,\varphi\,$ has no zeros. 
In this case
we can consider the conformally
changed metric $\tilde{g}:= \langle \varphi,\varphi \rangle ^{-2}g$.  Then the spinor field 
$\psi := |\langle
\varphi,\varphi \rangle|^{-1/2}\varphi$ is a twistor spinor of  constant length $\pm 1$ for 
$\tilde{g}$. Therefore,
we may
suppose that we have a twistor spinor $\varphi$ of constant  length $\pm 1$ on $(M,g)$.  
In that case $\langle
D\varphi,D\varphi\rangle$ is constant as well and $(M,g)$ is an Einstein space of scalar 
curvature   
\[ R= - \frac{4(n-1)}{n} \frac{\langle D\varphi,D\varphi \rangle}
{\langle \varphi,\varphi \rangle}\]
(cf. \cite{Baum/Friedrich/ua:91}).
If $R=0$ then either $\varphi$ or $D\varphi$ is a non-vanishing
parallel  spinor with timelike Dirac current. Hence, $(M,g)$ is a
product of the timelike line $\RR^{1,0}$ and a Riemannian spin manifold with  parallel spinor, 
which is then known to be
Ricci-flat and
K{\"a}hler. If $R$ is non-zero, either $\varphi$ is a Killing spinor or the spinors
\[ \psi_{\pm} = \frac{1}{2}\varphi \pm \sqrt{\frac{n-1}{nR}}D\varphi \]
are Killing spinors to different Killing numbers. In the first  case the
Killing number of $\varphi$ has to be
imaginary,  since the length of the spinor $\varphi$ is constant. Using Lemma \ref{lem-2}
one can show that (after normation and a change of
the time-orientation if the length of the spinor is +1) the  conditions of Proposition 
\ref{prop-4} are
satisfied and therefore, $(M,g)$ has to be Einstein--Sasaki.
In the second case when $\psi_{\pm}$ are imaginary Killing spinors
we can proceed in the same way. If $\psi_{\pm}$ are real  Killing spinors then $\psi_{\pm}$, 
$J\psi_{\pm}$ are four
$\CC$-linearly independent Killing spinors (here $J$  denotes the quaternionic structure 
on the spinor bundle) and,
by Lemma \ref{lem-3}, we can conclude that $(M,g)$ is conformally flat.\\
Let $n=7$ and let $\varphi$ be a twistor spinor without singularity. We consider the quaternionic structure
 $J:S \longrightarrow S$ on $S$ (Lemma \ref{lem-2}). In case
that $\langle\varphi,\varphi\rangle=0$ and $\langle\varphi,J\varphi\rangle=0$
the Dirac current $V_{\varphi}$ is lightlike and $(M,g)$ has to be conformally
equivalent to a Brinkmann space with parallel spinor. The restricted holonomy group of a 
7-dimensional Brinkmann space with parallel spinor is contained in $1\ltimes \RR^5 \subset SO(1,6)$,
$(SU(2)\times 1)\ltimes \RR^5 \subset SO(1,6))$ or $1\times SU(2) \subset SO(1,2)\times SO(4) \subset SO(1,6)$.
Hence, it is a pp-manifold, a Brinkmann space with special K\"ahler flag
or a metric product of a 3-dimensional
Minkowski space with a 4-dimensional Riemannian Ricci-flat K\"ahler manifold.\\
In case that
$\langle\varphi,\varphi\rangle\neq 0$ or
$\langle\varphi,J\varphi\rangle\neq 0$ at least one of the following
conformally changed metrics is
Einstein: $\tilde{g}=\langle\varphi,\varphi\rangle^{-2}g$,
$\tilde{g}=(Re\langle\varphi,J\varphi\rangle)^{-2}g$ or
$\tilde{g}=(Im\langle\varphi,J\varphi\rangle)^{-2}g$.
Therefore, we may assume that $(M,g)$ is an Einstein space if
the Dirac current $V_{\varphi}$ is not lightlike.
If $(M,g)$ is Ricci-flat, then $\varphi$ or $D\varphi$ are parallel with
timelike Dirac current.
Hence $(M,g)$ is the product of the timelike line $\RR^{1,0}$ and a
6-dimensional Riemannian Ricci-flat
K\"ahler manifold.
If the scalar curvature is non-zero, then $(M,g)$ admits a Killing spinor
$\varphi$.
First, we consider the case that the  
Killing number $\lambda$ is imaginary.
The spinor $\psi_{a,b}=a\varphi+bJ\varphi$ is an imaginary Killing
spinor with the same Killing number
for all $a,b\in \CC$. One easily checks that
$\,V_{\psi_{a,b}}=(|a|^2+|b|^2)\cdot
V_{\varphi}\,$. Therefore, the Clifford multiplication with $V_{\varphi}$
acts
as $\CC$-linear isomorphism on $span_{\CC}\{\varphi,J\varphi\}$.
This implies that there are complex numbers $\hat{a},\hat{b}$ such that
$V_{\psi_{\hat{a},\hat{b}}}\cdot\psi_{\hat{a},\hat{b}}=
\pm\psi_{\hat{a},\hat{b}}$. With Proposition \ref{prop-4} (after a
possible
change of time-orientation) we can conclude that $(M,g)$ is a Lorentzian
Einstein-Sasaki space.\\
Now, let $\varphi$ be a real Killing spinor on $(M,g)$ (i.e. the Killing
number is real) with timelike Dirac
current. Then the function $Q_{\varphi}= \langle \varphi,\varphi \rangle^2 + g(V_{\varphi},V_{\varphi})$ is 
constant on $M$. From Lemma \ref{lem-2} follows $\,Q_{\varphi}=-|\langle \varphi,J\varphi\rangle|^2 \leq 0\,$. 
In case $Q_{\varphi}=0$, $(M,g)$ is a warped product $(\RR\times F, -dt^2 + f(t)^2 h)$, where $(F,h)$ 
is a Riemannian manifold with parallel spinors and $f(t)=c \exp(\pm t)$ (cf. \cite{Bohle:03}). 
Hence, $(M,g)$ is locally conformal equivalent 
to a product of a timelike line with a 6-dimensional Riemannian Ricci-flat K\"ahler manifold. If $Q_{\varphi}<0$,
 then $(M,g)$ is a warped product $(\RR \times F, -dt^2 + f(t)^2 h)$, where $(F,h)$ is  
Riemannian manifold with real Killing spinors and $f(t)=d \cosh(t+c)$. Hence, $(M,g)$ is locally conformal equivalent to 
a product of a timelike line and a 6-dimensional nearly K\"ahler, non-K\"ahler manifold 
(cf. \cite{Baum/Friedrich/ua:91}). 
\qed
\end{proof}
\ \\
%


{\small

}

{\footnotesize Helga Baum, 
Institut f\"ur  Mathematik,
Humboldt-Universit\"at zu Berlin,
Sitz: Rudower Chaussee 25,
10099 Berlin,
Germany\\
email: baum@mathematik.hu-berlin.de\\
{\footnotesize Felipe Leitner,
University of Edinburgh, Department of Mathematics and Statistics, JCMB--King's Buildings, 
EH9 3JZ Edinburgh, Scotland}\\
email: felipe@maths.ed.ac.uk

\end{document}